\documentclass[10pt]{article}
\usepackage[latin1]{inputenc}
\usepackage{amsmath,amsthm,amssymb}
\usepackage{graphicx}
\usepackage{color}

\usepackage[hang]{caption2}
\usepackage{indentfirst}
\usepackage[section]{placeins}

\newtheorem{thm}{Theorem}
\numberwithin{thm}{section}
\newtheorem{lemma}[thm]{Lemma}

\newcommand{\neweq}[1]{\begin{equation}\label{#1}}

\def\phi{\varphi}
\def\RR{\mathbb R}
\def\CC{\mathbb C}
\def\NN{\mathbb N}

\def\ri{\rightarrow}

\def\incep{\left\{\begin{array}{ll} }
 \def\termin{\end{array}\right. }

\def\la2{\lambda^2}

 \newtheorem{theorem}{Theorem}[section]

\newcommand{\be}{\begin{equation}}
\newcommand{\ee}{\end{equation}}
\newcommand{\ba}{\begin{array}}
\newcommand{\ea}{\end{array}}
\newcommand{\bea}{\begin{eqnarray*}}
\newcommand{\eea}{\end{eqnarray*}}
\newcommand{\bean}{\begin{eqnarray}}
\newcommand{\eean}{\end{eqnarray}}

\newcommand{\ds}{\displaystyle}
\newcommand{\p}{\partial}
\newcommand{\s}{\sqrt}
\newcommand{\f}{\frac}

\newcommand{\ov}{\overline}
\newcommand{\q}{\quad}

\newcommand{\XX}{\phantom{XX}}

\begin{document}

\title{Existence of Frequency Modes  Coupling Seismic Waves and Vibrating Tall Buildings}
\author
{Darko Volkov and Sergey Zheltukhin \thanks{Department of Mathematical Sciences, Worcester Polytechnic
Institute,  Worcester MA 01609
({\tt darko@wpi.edu}, {\tt sergey@wpi.edu})} }

\maketitle

\begin{abstract}
We prove in this paper an existence result for 
frequency modes  coupling  seismic waves and vibrating tall buildings. 
 The derivation from physical principles of a set of equations  modeling this phenomenon
was done in previous studies. In this model all vibrations are assumed to be anti plane and time harmonic so the two dimensional Helmholtz equation
can be  used.  
A coupling frequency mode is obtained once we can
 determine a wavenumber such that the solution of the corresponding Helmholtz equation
in the lower half plane with relevant Neumann and Dirichlet at the interface satisfies a specific integral
equation at the base of an idealized tall building. 
Although  numerical simulations suggest 
that such wavenumbers should exist, as far as we know, to date, there is no theoretical proof of existence. 
This is what  this present study offers to provide.

\end{abstract}

\normalsize

\section{Introduction}
The traditional approach to evaluating seismic risk in urban areas is to consider  seismic waves  under ground 
 as the only cause for  motion above ground. In earlier  studies, seismic wave propagation was evaluated 
in an initial step and in a second step  impacts on  man made structures  were inferred. However, observational 
evidence has since then suggested that when an earthquake strikes a large city,  seismic activity 
may in turn be
 altered by the  response of  the buildings. 
This phenomenon is referred to as the ``city-effect''
and has been studied by many authors, see \cite{shuttles,ullevi}.\\
More recently, \cite{ghergu},
 Ghergu and Ionescu
have derived a model for the city effect based on the equations of solid mechanics and appropriate
coupling of the different elements involved in the physical set up of the problem.
They then proposed a clever way to compute a numerical solution to their system of equations.
In this way, in \cite{ghergu},
Ghergu and Ionescu were able to compute a city frequency constant:
 given the geometry and the specific physical constants of an idealized two dimensional city, 
they computed a frequency that  leads to  coupling between vibrating buildings and underground seismic waves. 
In this present paper our goal is to prove that 
the equations modeling the city effect introduced in \cite{ghergu} 
are solvable. 
We acknowledge that these equations were carefully derived following the laws of solid mechanics
combined to the knowledge of  the relevant dominant effects causing this phenomenon. There is also
ample numerical evidence that these coupling frequencies should exist, see
\cite{ghergu,DVandSergey_computation}, at least in the range of physical parameters under consideration
in these numerical simulations. 
As far as we know, there is no mathematical proof, however, that these coupling frequencies must exist.
To fill this gap, we will show in this paper that given a building with positive height, mass, and elastic modulus,
the (rather involved) set of coupled equations given in \cite{ghergu} 
defining frequencies
coupling that building and the ground beneath have at least one solution, if some constants
coming from non-dimensionalization of physical parameters satisfy a sign condition. 
Furthermore, we show that once
the constants for the physical properties of the ground and
the building are fixed,    the set of all possible coupling frequencies is finite.\\
Here is an outline of this paper. In
section 2 we introduce the equations defining frequency modes coupling seismic waves and vibrating tall 
buildings. For the sake of brevity we directly provide them in non-dimensional form. A derivation from
physical principles and non-dimensionalization calculations can be found in \cite{ghergu,DVandSergey_computation}.
In section 2 we introduce on the one hand the PDE modeling the propagation of time harmonic
waves under ground while at the ground level  a building is subject to a given displacement 
and a no force condition is applied elsewhere, and on the other hand the integral equation
ensuring coupling between vibrations under ground and in the building.
 The set of coupled equations for which we prove existence is comprised of the PDE and of the coupling integral 
equation. \\
In section 3 we have to carefully study the boundary Dirichlet to Neumann operator $T_k$ for Helmholtz problems outside the 
unit disk in $\RR^2$, where $k$ is the wavenumber. We are aware that this is a well known operator, however,
for our purposes, we need to show the lesser known fact that that $T_k$ is (strongly) real analytic in $k$, and we need
to determine the strong limit of $T_k$ as $k$ tends to zero. 
In section 4 we reformulate our half plane problem to the whole plane using a symmetry: that way the operator
$T_k$ introduced in section 3 can be used. Since the strong limit of $T_k$ was found in  section 3,
it is then possible to study the low frequency behavior of our problem thanks to 
manipulations of $k$- dependent variational problems. \\
Section 5 deals with the much more delicate question of high frequency asymptotics.  Understanding how waves
behave at high frequencies has always piqued the interest of scientists. The geometric optics approximation has been known
for quite some time; in the late 19th century Kirchhoff wrote down specific equations
capturing the behavior of waves at high frequencies. A rigorous mathematical study of these phenomena first appeared
in papers by Majda, Melrose and Taylor, see 
\cite{Majda,MajdaTaylor,MelroseTaylor}. We note, however, that their results are limited to the case 
where scatterers are convex domains, 
so their results not applicable to our particular case.
More recently, Chandler-Wildez, Hewett,  and, Langdony, see \cite{Chandler1,Chandler2},
published continuity and coercivity estimates 
 pertaining to  either
 scattering in dimension 2 by soft or hard line segments (our case), or
scattering in dimension 3 by soft or hard  open planar surfaces.
These estimates include bounds whose 
\textbf{explicit dependence on wavenumbers} is stated and proved.
In section 5, after  first  informally deriving  the expected behavior of some quantities relevant to 
the coupling frequency problem, 
we turn the rigorous proof of that expected result. This is where
the new estimates by Chandler et al. turn out to be crucial.
 Finally, in section 6, we combine all the intermediate results obtained in  previous sections 
to complete the proof of our main theorem. 
In section 7, there ensues a brief discussion on our findings and on how we plan to extend 
this present study to more complex geometries  in future work. 
This paper also contains an appendix with  an overview of results on Hankel functions relevant to our
work.


\section{The equations defining frequency modes coupling seismic waves to vibrating tall 
buildings. Statement of main theorem.}
Following \cite{ghergu} and \cite{DVandSergey_computation}  we model the ground 
to be 
the  elastic half-space 
$x_2>0$ in three dimensional  space, where $(x_1, x_2, x_3)$ is the space variable.
We only considered the \textbf{anti-plane shearing} case:
all displacements  occur in the $x_3$ direction 
and are independent of $x_3$.
Since in the rest of this paper we won't be using the third direction,
we set $x=(x_1, x_2)$.
We denote by $\RR^{2+}$ the half plane $x_2>0, \, x_3=0$.
We refer to \cite{ghergu} and \cite{DVandSergey_computation} for
a careful derivation of how given the mass density, the shear rigidity
of the building and the ground, the height and the width of the building,
the mass at the top and the mass at the foundation of the building, after non dimensionalization,
we arrive at the following  system of equations, 
assuming that the building has  rescaled width $1$  
and is standing on the
$x_1$ axis, so that its foundation may be assumed to be the line segment
$\Gamma = [-\f12,\f12] \times \{0\}$,
\bean
  \label{non_dim_helmholtz}
  \Delta \Phi + k^2 \Phi = 0 \mbox{ in } \RR^{2+}, \\
  \Phi=1 \mbox{ on } \Gamma,
  \label{non_dim_boundary} \\
  \displaystyle \frac{\partial \Phi}{\partial x_2} = 0 \mbox{ on } \{ x_2=0 \} \setminus \Gamma,\\
  \f{\p \Phi}{\p r} - ik \Phi = o(r^{-1}), \mbox{ uniformly as } r \ri \infty
\label{uext0}
  \\ \label{non_dim_eigenvalue}
  \displaystyle q(k^2)  = p(k^2) \mbox{Re }\int_{\Gamma} {\frac{\partial \Phi} {\partial x_2} (s,0)} ds 
\eean
where
\begin{equation}
\label{p_and_q}
  p(t) = C_1 t -C_2, \quad q(t)= t (C_3 t + C_4)
\end{equation}
Here $k>0$ is the wavenumber, 
$r= \sqrt{x_1^2+x_2^2}$, 
the rescaled physical displacement is $\mbox{Re } \Phi e^{-i k t}$, 
and the  constants $C_1, C_2, C_3, C_4$, are determined by the physical 
properties of the underground and the building as specified in 
\cite{ghergu,DVandSergey_computation}.
Note that system (\ref{non_dim_helmholtz}-\ref{p_and_q}) is 
\textbf{non linear} in the unknown wavenumber $k$.
The goal of this paper is to show the following theorem,

\begin{thm} \label{main_th}
For any positive value of  the  constants $C_1, C_2, C_3, C_4$, the system of equations
(\ref{non_dim_helmholtz}-\ref{p_and_q}) has at least one solution in $k$, that 
is, there exists a positive $k$ and a function $\Phi$ which has 
locally $H^1$ regularity in $\ov{\RR^{2+}}$ such that equations
(\ref{non_dim_helmholtz}-\ref{p_and_q}) are satisfied.
Moreover, this equation has at most a finite number of solutions.
\end{thm}

Standard arguments can show that if we fix a positive $k$ the system of equations
(\ref{non_dim_helmholtz}) through (\ref{uext0}) is uniquely solvable. 
Theorem \ref{main_th} asserts that 
for some of those $k$'s the additional relation (\ref{non_dim_eigenvalue}) will hold.
Here is a sketch of the proof of theorem \ref{main_th}.
For  $k$ in $(0, \infty)$ we define
\bean
\label{Fdef}
F(k) =q(k^2)  - p(k^2) \mbox{Re }\int_{\Gamma} {\frac{\partial \Phi_{k}} {\partial x_2} (s,0)} ds 
\eean 
where $\Phi_k$ solves (\ref{non_dim_helmholtz}) through (\ref{uext0}).
We will first show that $F$ is real analytic in $k$. 
Then we will perform a low frequency and a high frequency  analysis of $\Phi_k$.
The low frequency analysis will show that $F$ must be negative 
in $(0, \alpha)$ for some positive number $\alpha$.
The high frequency analysis will prove that $\ds \lim_{k \ri \infty} F(k) =  \infty$,
concluding the proof of theorem \ref{main_th}.

\section{The boundary Dirichlet to Neumann operator: analicity with regard to the wavenumber}

We denote by $D$ the open unit disk of $\RR^2$ centered at the origin. 
Our first lemma is certainly well known to most readers, but we chose
to include it since it is helpful in the detailed study of the related
Dirichlet to Neumann operator relevant to this study.
\begin{lemma}\label{extend}
Let $k>0$ be a wave number and $f$ be a function in
the Sobolev space $\ds H^{\f12}(\p D)$.
The problem
\bean
(\Delta + k^2) u =0 \mbox{ in } \RR^2 \setminus \ov{D} \label{uext1}\\
u=f \mbox{ on } \p D \label{uext2} \\
\f{\p u}{\p r} - ik u = o(r^{-1}), \mbox{ uniformly as } r \ri \infty
\label{uext3}
\eean
has a unique solution.
Writing $f = \ds \sum_{n=-\infty}^{\infty} a_n e^{in \theta}$, we have 
\bean
u = \ds \sum_{n=-\infty}^{\infty} a_n e^{in \theta} 
\f{H_n (k r)}{H_n(k)},
\label{useries}
\eean
  This series and all its derivatives are uniformly convergent on any  subset of $\RR^2$ in the form $r \geq A $ where $A>1$.
\end{lemma}
\textbf{Proof:}
Existence and uniqueness for equation (\ref{uext1}-\ref{uext3}) are well known, we are chiefly interested here in the convergence properties of the series (\ref{useries}). We first note that $H_{-n}(z)=(-1)^n H_n(z)$, so we will establish convergence properties as $n \ri \infty$ and the case $n \ri -\infty$ will then follow easily. From formula (\ref{largen})
in appendix,  we can claim that
$$
\f{H_n (k r)}{H_n(k)} \sim \f{1}{r^n},
$$ 
uniformly in $r$ as long as $r$ remains in a compact set of $(0,\infty)$. Fix two real numbers $A, B$ such that $1<A<B$. Set $M = \sup |a_n|$. It follows that 
$$
| a_n e^{in \theta} 
\f{H_n (k r)}{H_n(k)}| \leq 2 \f{M}{A^n}
$$
for all $n$ large enough, uniformly for all $r$ in $[A,B]$, so the series (\ref{useries}) is uniformly convergent on any compact subset of $\RR^2 \setminus \ov{D}$.\\
Next we use the recurrence formula (\ref{besselrec}) given in appendix to write
$$
k\f{H_n' (k r)}{H_n(k)} =- k\f{H_{n+1} (k r)}{H_n(k)}
+ \f{n}{ r } \f{H_n (k r)}{H_n(k)}
$$
The term $\ds \f{n}{ r } \f{H_n (k r)}{H_n(k)}$ can be estimated as previously. For $\ds k\f{H_{n+1} (k r)}{H_n(k)}$ use 
 (\ref{largen}) one more time to find
$$
\f{H_{n+1} (k r)}{H_n(k)} \sim \f{2n+2}{r^{n+2}}
$$ 
uniformly for all $r$ in $[A,B]$. 
At this stage we can conclude that the $r$ derivative of the series (\ref{useries}) is uniformly convergent on any compact subset of $\RR^2 \setminus \ov{D}$.\\
A $\theta$ derivative of the series (\ref{useries}) corresponds to a multiplication by $ in$ so the uniform convergence property holds for that derivative too. Second derivatives can be treated in a similar way to find that the function defined by the series (\ref{useries}) is in $C^2(\RR^2 \setminus \ov{D})$. \\
If $r \geq 2$, we combine lemma \ref{decrease} and (\ref{largen}) (in appendix) to write
$$
|\f{H_n (k r)}{H_n(k)}| \leq |\f{H_n (2k )}{H_n(k)}|
\leq 2^{-n+1},
$$
for all $n$ greater than some $N$, for all $r \geq 2$. Similarly
$$
|k\f{H_n' (k r)}{H_n(k)} | \leq | k\f{H_{n+1} (2k )}{H_n(k)}|
+ |\f{n}{ r } \f{H_n (2k )}{H_n(k)}| \leq n 2^{-n+1},
$$
for all $n$ greater than some $N$, for all $r \geq 2$. Given that $a_n$ is bounded, it follows that  the series (\ref{useries}) and its $r$ derivative are uniformly convergent for all $r$ 
in $[2, \infty)$. A similar argument can be carried out for the $\theta$ derivative, and all second derivatives. \\
Next we recall that $H_n$ satisfies the Bessel differential equation
$$
 y''(r) + \f{1}{r} y'(r) + (1 - \f{n}{r^2}) y(r) = 0 
$$
to argue that each function $e^{i n \theta} H_n(k r)$ 
satisfies Helmholtz equation due to the form of the Laplacian in polar 
coordinates, namely $\ds \p_r^2 + \f{1}{r} \p_r + \f{1}{r^2} \p_\theta^2$.
All together this shows that  the function defined by the series
(\ref{useries})  satisfies (\ref{uext1}). \\
To prove (\ref{uext3}) we first note that each function $e^{i n \theta} H_n(k r)$ satisfies that estimate due to the well known asymptotic behavior of Hankel's functions $H_n$, see \cite{Abra}. From there (\ref{uext3}) can be derived using that the series (\ref{useries})
and its $r$ derivative are uniformly convergent for all $r$ in $[2, \infty)$.  \\
Finally it is worth mentioning that for any fixed $r \geq 1$ the series in (\ref{useries}) is in the Sobolev space $H^{\f12}(\p D)$ for all $r \geq 1$ and that further applications of lemma \ref{decrease} will show that this series converges strongly to $f$ in $H^{\f12}(\p D)$  as $r \ri 1^+$. $\square$

We now define the linear operator $T_k$ which maps $H^{\f12}(\p D)$ into $H^{-\f12}(\p D)$ by the formula
\bean
T_k (f) = \ds \sum_{n=-\infty}^{\infty} a_n e^{in \theta} 
k \f{H_n' (k )}{H_n(k)},
\eean
where $\ds f =\sum_{n=-\infty}^{\infty} a_n e^{in \theta} $.
$T_k$ is continuous since formula  (\ref{besselrec}) combined to (\ref{largen}) 
(both given in appendix) implies that
$$
k \f{H_n' (k )}{H_n(k)} \sim -n , \q n\ri \infty
$$
According to lemma \ref{extend},
an equivalent way of defining $T_k$ is to say that it maps $f$ to $\ds \f{\p u }{\p r}|_{r=1}$, where $u$ is the solution to
(\ref{uext1} - \ref{uext3}).
We denote by $<,>$ the duality bracket between
$H^{\f12}(\p D)$  
and $H^{-\f12}(\p D)$ which extends the dot product 
$<f,g>=\int_{\p D} f \ov{g}$.
\begin{lemma}\label{signTk}
 Let $f$ be  in $H^{\f12}(\p D)$.
Then $\mbox{Re} <T_k(f), f> \leq 0$. 
\end{lemma}
\textbf{Proof:} \\
Set $\ds f =\sum_{n=-\infty}^{\infty} a_n e^{in \theta} $. By definition of $T_k$
$$
\mbox{Re} <T_k(f), f> = 2 \pi \sum_{n=-\infty}^{\infty} |a_n|^2
k \mbox{Re}( \f{H_n' (k )}{H_n(k)})
= 2 \pi \sum_{n=-\infty}^{\infty} |a_n|^2
k  \f{\mbox{Re}(H_n' (k ) \ov{H_n(k)})}{|H_n(k)|^2}
$$
The result then follows from lemma \ref{decrease} given in appendix.
$\square$

\begin{lemma} \label{Tkanal}
$T_k$ is real analytic in $k$ for
$k$ in $(0,\infty)$.
\end{lemma}
\textbf{Proof:}
Let $k_0$ and $b$ be two real numbers such that $0<b<k_0$. Define $D_b(k_0)$ the closed disk in the complex plane centered at $k_0$ and with radius $b$. According to \cite{Kato}, chapter 7,  to show that $T_k$ is real analytic in $k$ for $k$ in $(0,\infty)$ in operator norm, it suffices to
fix $f$ and $g$ in  $H^{\f12}(\p D)$ and to show that $<T_k(f),g>$ is an analytic function of $k$ in $D_b(k_0)$. Writing  $\ds f =\sum_{n=-\infty}^{\infty} a_n e^{in \theta} $, $\ds g =\sum_{n=-\infty}^{\infty} b_n e^{in \theta} $, 
we have
\bean \label{sumforanal}
<T_k(f),g> = 2\pi \sum_{n=-\infty}^{\infty} a_n \ov{b_n}
k \f{H_n' (k )}{H_n(k)}
\eean
Given that $\ds \sum_{n=-\infty}^{\infty} n|a_n b_n| <\infty$ and
$\ds k \f{H_n' (k )}{H_n(k)} \sim -|n| , \q |n|\ri \infty$, uniformly for all $k$ in  $D_b(k_0)$
(see lemma \ref{realeq2} in appendix), 
the series in (\ref{sumforanal}) is a uniformly convergent sum of analytic functions of $k$, thus $<T_k(f),g>$ is analytic in  $D_b(k_0)$.
$\square$

\subsection{The limit of the boundary operator $T_k$ as $k$ approaches zero}
\begin{lemma}\label{Tconv}
The operator $T_k$ converges strongly to the operator $T_0$ which maps $H^{\f12}(\p D)$  
into $H^{-\f12}(\p D)$ and is defined by the formula
\bean 
  \label{T0def}
T_0 (f) = \ds \sum_{n=-\infty}^{\infty}-|n| a_n e^{in \theta},
\eean
where $\ds f =\sum_{n=-\infty}^{\infty} a_n e^{in \theta} $.
\end{lemma}
\textbf{Proof:}
Let  $\ds f =\sum_{n=-\infty}^{\infty} a_n e^{in \theta} $ be in 
 $H^{\f12}(\p D)$ .
We can write
$$
\| T_k(f) - T_0(f) \|_{H^{-\f12}(\p D)} =
\sum_{n=-\infty}^{\infty} \f{|a_n|^2}{\s{n^2+1}}
|k \f{H_n' (k )}{H_n(k)} - |n||^2
$$
and then apply lemma \ref{realeq2} given in appendix. $\square$

\begin{lemma}\label{laplace}
Let  $f$ be a function in the Sobolev space $H^{\f12}(\p D)$. The problem
\bean
\Delta u =0 \mbox{ in } \RR^2 \setminus \ov{D} \label{uext1lap}\\
u=f \mbox{ on } \p D \label{uext2lap} \\
 u = O(1), \mbox{ uniformly as } r \ri \infty
\label{uext3lap}
\eean
has a unique solution. Writing $f = \ds \sum_{n=-\infty}^{\infty} a_n e^{in \theta}$,
we have 
\bean
u = \ds \sum_{n=-\infty}^{\infty} a_n e^{in \theta} r^{-|n|},
\label{userieslap}
\eean
 This series and all its derivatives are uniformly convergent
on any  subset of $\RR^2$ in the form $r \geq u $ where $u>1$.
\end{lemma}
\textbf{Proof:}
Existence and uniqueness for equation (\ref{uext1lap}-\ref{uext3lap}) are well known, see \cite{Folland}. Proving the uniform convergence properties is trivial given how simple the radial terms are. $\square$\\\\
We note that an equivalent way of defining $T_0$ is to say that it maps $f$ to $\ds \f{\p u }{\p r}|_{r=1}$ where $u$ is the solution to
(\ref{uext1lap}-\ref{uext3lap}).

\section{An equivalent problem in the plane $\RR^2$ minus a line segment. Low wavenumber asymptotics.}
\label{section_equivalent_problem}
We use the following notation in this section: $\Gamma$ is the line segment $\ds [-\f12, \f12] \times \{0\}$, $\Omega$ is the open set $\{x \in \RR^2: |x| <1\} \setminus \Gamma$. $\p D$ will denote the same boundary as in the previous section. On $\Gamma$ an upper and a lower trace for all functions
in  $H^1(\Omega)$ can be defined:  even though $\Gamma$ is an interior boundary of $\Omega$, since $\Gamma$ is $C^1$, an adequate version of the trace theorem on $\Gamma$ holds, see
\cite{Troianiello}. \\
We denote by $H^1_{0,\Gamma}(\Omega)$ the closed subspace of $H^1(\Omega)$ consisting of functions whose upper and lower trace on $\Gamma$ are zero.
\begin{lemma}
Let $L$ be a continuous linear functional on $H^1_{0,\Gamma}(\Omega)$. The following variational problem has a unique solution for any $k>0$:\\
find $u$ in  $H^1_{0,\Gamma}(\Omega)$ such that
\bean
\label{var}
\int_\Omega \nabla u \cdot \nabla v  - k^2 uv - \int_{\p D} (T_k u)\, v = L(v),
\eean
for all $v$ in $H^1_{0,\Gamma}(\Omega)$.
\end{lemma}
\textbf{Proof:}
We first prove uniqueness. Assume that $u$ is in $\ds H^1_{0,\Gamma}(\Omega)$ and satisfies (\ref{var}) with $L=0$. It is clear that $u$ satisfies $(\Delta + k^2) u =0$ in $\Omega$. Define
$$
\Omega^+ = \{ (x_1,x_2) \in \Omega: x_2 >0 \}
$$
and define $\Omega^-$ likewise. By Green's theorem we must have that 
$$
\mbox{Im}\int_{\p \Omega^+} u \ov{\f{\p u}{\p n}}=
\mbox{Im}\int_{\p \Omega^-} u \ov{\f{\p u}{\p n}} =0
$$
Using the fact that the upper and lower traces of $u$ on $\Gamma$ are zero we infer that
\bean \label{forRei}
\mbox{Im}\int_{\p D} u \ov{\f{\p u}{\p n}} =0.
\eean
Next we observe that the variational problem (\ref{var}) implies
\bean \label{pDcont}
\f{\p u}{\p n}= T_k u
\eean
on the boundary $\p D$. Outside $D$ we extend $u$ by setting it equal to the solution
of (\ref{uext1}-\ref{uext3}) where $f$ is the trace of $\ds u_{| \Omega}$ on $\p D$. We can claim thanks to (\ref{pDcont}) that $u$ and $\ds \f{\p u}{\p n}$ are continuous across $\p D$. By (\ref{forRei}) $u $ must be zero in $\RR^2\setminus D$ due to Reillich's lemma. Consequently, $u$ and $\ds \f{\p u}{\p n}$ are zero on $\p D$, and therefore $u$ is also zero in $\Omega$ due to the Cauchy Kowalewski theorem. \\ 
To show existence define the bilinear functional
\bean \label{akdef}
a_k(u,v) = \int_\Omega \nabla u \cdot \nabla v  - k^2 uv -
\int_{\p D} (T_k u)\, v
\eean
$a_k$ is continuous on $H^1_{0,\Gamma}(\Omega) \times
H^1_{0,\Gamma}(\Omega)$. Due to lemma \ref{signTk}
we can claim that
$$
\mbox{Re} \, a_k(u,\ov{u}) \geq \|\nabla u \|^2_{L^2(\Omega)} 
- k^2 \|u \|^2_{L^2(\Omega)}
$$
We  note that thanks to a generalized version of Poincare\'e's
inequality
$$
v \ri \|\nabla v \|_{L^2(\Omega)}
$$
is a norm on  $H^1_{0,\Gamma}(\Omega)$ which is equivalent to the
natural norm.
Now, since the injection of $H^1_{0,\Gamma}(\Omega)$
into $L^2(\Omega)$ is compact we can claim that either the equation $a_k(u,v) = L(v)$ is  uniquely solvable or the equation $a_k(u,v)=0$ has non trivial solutions. Given that we proved uniqueness, we conclude that $a_k(u,v) = L(v)$ is uniquely solvable and that $u$ depends continuously on $L$. $\square$

\begin{lemma}
Let $L$ be a continuous linear functional on $H^1_{0,\Gamma}(\Omega)$.
The following variational problem has a unique solution:\\
find $u$ in  $H^1_{0,\Gamma}(\Omega)$ such that
\bean\label{var0}
\int_\Omega \nabla u \cdot \nabla v  -  
\int_{\p D} (T_0 u)\, v = L(v),
\eean
for all $v$ in $H^1_{0,\Gamma}(\Omega)$.
\end{lemma}
\textbf{Proof:}
We observe that due to the definition (\ref{T0def}) of $T_0$,
$<T_0(f), f> $ is real for all $f$ in $H^{\f12} (\p D)$ 
and $<T_0(f), f> \leq 0$. 
We conclude that problem (\ref{var0})
is uniquely solvable and that the solution $u$ depends continuously
on $L$. $\square$

Let $\phi$ be a smooth compactly supported function in $D$
which is equal to 1 on $\Gamma$ and such that $\phi(x_1, -x_2)=
\phi(x_1,x_2)$.
For all $k \geq 0$ we set $\tilde{u}_k$ in $H^1_{0,\Gamma}(\Omega)$
 to be the solution to
\bean\label{varphi}
\int_\Omega \nabla\tilde{u}_k \cdot \nabla v  - k^2 \tilde{u}_k v -
\int_{\p D} (T_k \tilde{u}_k)\, v = \int_\Omega (\Delta \phi + k^2\phi)v,
\q \forall  v \in H^1_{0,\Gamma}(\Omega)
\eean
and we set $u_k = \tilde{u}_k + \phi$. 

\begin{lemma} \label{onuk}
For $k>0$, $u_k $
satisfies the following properties:\\
(i). $u_k$ is in $H^1(\Omega)$. \\
(ii). The upper and lower trace on $\Gamma$ of $u_k$
are both equal to the constant 1.\\
(iii). $u_k$ can be extended to a function in $\RR^2 \setminus\Gamma$
such that, if we still denote by $u_k$ that extension, \\
\XX $(\Delta + k^2) u_k =0$ in   $\RR^2 \setminus\Gamma$, \\ 
\XX $\ds \f{\p u_k}{\p r} - ik u_k = o(r^{-1}), \mbox{ uniformly as } r \ri \infty$ \\
\XX $u_k(x_1,-x_2) =u_k(x_1,x_2) $ for all $(x_1,x_2)$ in $\RR^2 \setminus\Gamma$\\
\XX $\ds \f{\p u_k}{\p x_2} (x_1 ,0) =0$ if $(0,x_1) \notin \Gamma$.\\
(iv). Denoting $\ds \f{\p u_k}{\p x_2^-}$ the lower trace of
 $\ds \f{\p u_k}{\p x_2}$ on $\Gamma$,
\bean \label{intnon0}
 \mbox{Im}  \int_\Gamma  \f{\p u_k}{\p x_2^-} <0
\eean 
\end{lemma}
\textbf{Proof:}
Properties (i) and (ii) are clear. The first two items of property (iii) hold simply because we can write $u_{| \p D}= \ds \sum_{n=-\infty}^{\infty} a_n e^{in \theta}$ and then set $u = \ds \sum_{n=-\infty}^{\infty} a_n e^{in \theta} \f{H_n (k r)}{H_n(k)}$ in $\RR^2 \setminus \ov{\Omega}$. We then use Lemma \ref{extend} in combination to the fact that variational problem (\ref{varphi}) implies that $T_k u_k$ is the limit of $\ds \f{\p u_k}{\p r}$ as $r \ri 1^-$. \\
To show the third item in (iii) we set $\ds \underline{\tilde{u_k}} (x_1, x_2) = \tilde{u_k} (x_1, -x_2) $ and for any  arbitrary $v$ in $H^1_{0,\Gamma}(\Omega)$, $\ds \underline{v} (x_1, x_2) =v (x_1, -x_2) $. Next we observe that 
\bea
\int_\Omega \nabla\underline{\tilde{u_k}} \cdot \nabla v  - k^2 \underline{\tilde{u_k}} v  = 
\int_\Omega \nabla\tilde{u}_k \cdot \nabla \underline{v}  - k^2 \tilde{u}_k \underline{v}
\eea
In polar coordinates we have the relations $ \ds \underline{\tilde{u_k}} (r, \theta) = \tilde{u_k} (r, -\theta) $ and $\ds \underline{v} (r, \theta) =v (r, -\theta) $, so 
\bea
\int_{\p D} (T_k \underline{\tilde{u_k}})\, v = \int_{\p D} (T_k 
\tilde{u_k})\, \underline{v} 
\eea
Finally, since $\phi$ is even in $x_2$,
\bea
\int_\Omega (\Delta \phi + k^2\phi) \underline{v} = 
\int_\Omega (\Delta \phi + k^2\phi)v
\eea
Since the solution to problem (\ref{varphi}) is unique we must have
$\underline{\tilde{u_k}} = \tilde{u_k}$, proving the third item in
(iii). \\
Since $u_k$ is even in $x_2$, it follows that $\ds \f{\p u_k}{\p x_2}$ is zero on the line $x_2 =0$ minus the segment $\Gamma$, proving the last item in (iii).\\
Since $u=1$ on $\Gamma$,
$$
\mbox{Im} \int_\Gamma  \f{\p u_k}{\p x_2^-} =
\mbox{Im} \int_\Gamma  \f{\p u_k}{\p x_2^-} \ov{u_k} =
-\mbox{Im} \int_{\p D^-}  \f{\p u_k}{\p r} \ov{u_k} 
$$
where $\p D^-$ is the intersection of the circle $\p D$ and the lower
half plane $x_2<0$.
We use parity one more time to argue that 
$$\mbox{Im} \int_{\p D^-}  \f{\p u_k}{\p r} \ov{u_k} =
\f12  \mbox{Im} \int_{\p D}  \f{\p u_k}{\p r} \ov{u_k}
$$
Knowing that $u_k$ is not zero everywhere and 
using Reillich's lemma 
we can claim that $\ds \mbox{Im} \int_{\p D} \f{\p u_k}{\p r} \ov{u_k} >0$. $\square$

\begin{lemma}\label{convth}
Let $\tilde{u}_k$ be the solution to (\ref{varphi}) for all $k\geq0$,
and set $u_k = \tilde{u}_k + \phi$ . \\
(i). $u_k$ is analytic in $k$ for $k>0$. \\
(ii). $u_k$ converges strongly to $u_0$ in $H^1_{0,\Gamma}(\Omega)$. More precisely there is a constant $C$ such that 
\bean \label{ukconvest}
 \|u_k- u_{0} \|_{H^1(\Omega)} \leq C (k^2 + \| T_k - T_0 \|)
\eean
\end{lemma}
\textbf{Proof:}
To show (i) we define an operator $A_k$
from $H^1_{0,\Gamma}(\Omega)$ to its dual defined by
$$
<(A_k u), v> = a_k(u,v), \q \forall v \in H^1_{0,\Gamma}(\Omega)
$$
where $a_k$ was defined in (\ref{akdef}).
Thanks to lemma \ref{Tkanal} we can claim that $A_k$ is analytic
in $k$ for $k>0$. But we showed $A_k$ is invertible for $k>0$ and that
its inverse is a continuous linear functional.  According to \cite{Kato}, chapter 7, $A_k^{-1}$ is then also analytic
in $k$ for $k>0$, and so is $A_k^{-1}L$ for any fixed $L$
 in the dual of $H^1_{0,\Gamma}(\Omega)$.
Note that we have obtained analicity of $u_k$ relative to the
$H^1(\Omega)$ norm.  \\
To prove statement (ii)
we first show that $\|\tilde{u}_k \|_{H^1(\Omega)} $
is uniformly bounded for all $k$ in $[0,B]$ where $B$ is a positive 
constant. 
Set $v=\tilde{u}_k $
in (\ref{varphi}) and use lemma \ref{signTk} to obtain
\bean \label{unibound}
\|\nabla  \tilde{u}_k \|^2_{L^2(\Omega)} \leq &
k^2 \| \tilde{u}_k \|^2_{L^2(\Omega)}
+ C \| \phi \|_{H^2(\Omega)}  \|\tilde{u}_k \|_{L^2(\Omega)},  
\eean
where $C$ is a positive constant.
We now need to invoke Poincare\'e 's inequality which in $H^1_{0,\Gamma}(\Omega)$
implies that there is a positive constant $C_p$ such that
\bean \label{poincarre}
\| v \|_{L^2(\Omega)} \leq C_p \| \nabla v \|_{L^2(\Omega)},
\eean
for all $v$ in $H^1_{0,\Gamma}(\Omega)$.
Since (\ref{unibound}) implies that
\bean \label{unibound2}
\|\nabla  \tilde{u}_k \|^2_{L^2(\Omega)} \leq 
\f{k^2 C_{p}^2}{2} + \f{1}{2 C_p^2}\| \tilde{u}_k \|^2_{L^2(\Omega)}
+ C^2 C_p^2 \| \phi \|_{H^2(\Omega)}^2+ \f{1}{4 C_p^2}  \|\tilde{u}_k \|_{L^2(\Omega)},  
\eean
we may now use (\ref{poincarre}) to conclude that
   $\|\tilde{u}_k \|_{H^1(\Omega)} $ is bounded 
for all $k$ in $[0,B]$, as long as $B$ is small enough.\\
Next
we note that $u_k - u_0 = \tilde{u}_k - \tilde{u}_0 $
and satisfies for all $v$ in $H^1_{0,\Gamma}(\Omega)$
\bea
\int_\Omega \nabla (u_k  - u_0) \cdot \nabla v  -\int_\Omega k^2 \tilde{u}_k v -
\int_{\p D} (T_k \tilde{u}_k - T_0 \tilde{u}_0)\, v = 
\int_\Omega k^2 \phi v
\eea
We re write $\int_{\p D} (T_k \tilde{u}_k - T_0 \tilde{u}_0)\, v$ 
as
\bea \ds
\int_{\p D}\ds  \left(T_k (\tilde{u}_k- \tilde{u}_{0}) \right) v -
 \int_{\p D} \left( (T_0 - T_k ) \tilde{u}_0 \right) v\, 
\eea
we choose $v =\ov{\tilde{u}_k- \tilde{u}_{0}} $ 
and we use that, due to lemma \ref{signTk},
$$
\mbox{Re} \int_{\p D}\ds  \left(T_k (\tilde{u}_k- \tilde{u}_{0}) \right) v
\leq 0
$$
to infer the inequality
\bea
\|\nabla( u_k- u_{0}) \|^2_{L^2(\Omega)} \leq &
k^2 \| \tilde{u}_k \|_{L^2(\Omega)}
\| u_k- u_{0} \|_{L^2(\Omega)}  \\&+ 
\| T_k - T_0 \| \| \tilde{u}_0 \|_{H^{\f12} (\p D)}
\| u_k- u_{0} \|_{H^{\f12} (\p D)} 
+ k^2 \| \phi \|_{L^2(\Omega)}  \|u_k- u_{0} \|_{L^2(\Omega)}  
\eea
the result follows since we know that $ \| \tilde{u}_k \|_{L^2(\Omega)}$ is bounded for $k$ in $[0,B]$ and by application 
of the trace theorem and of Poincare\'e 's inequality (\ref{poincarre}). $\square$.

\begin{lemma}
Denote $\ds \f{\p u_k}{\p x_2^\pm}$ the upper and lower traces of
$\ds \f{\p u_k}{\p x_2}$ on $\Gamma$. Then \\
(i). 
\bean \label{opp}
\ds \f{\p u_k}{\p x_2^+} = - \f{\p u_k}{\p x_2^-}
\eean
(ii). Denote $G_k(x,y) = \f{i}{4} H_0 (k |x-y|)$. For all $x$ in $\Omega$,
\bean \label{intform}
u_k (x) = 2 \int_\Gamma  G_k(x,y) \f{\p u_k}{\p x_2^-} (y) dy.
\eean
\end{lemma}
\textbf{Proof:}
Denote $\Omega^+ = \{(x_1,x_2) \in \Omega: x_2 >0 \}$ and $\Omega^- = \{(x_1,x_2) \in \Omega: x_2 <0 \}$. It is well known from potential theory that if $x$ is in $\Omega^+$
\bean
u_k (x) =  \int_{\p \Omega^+} G_k(x,y) \f{\p u_k}{\p n} (y) 
- \f{\p  G_k(x,y)}{\p n(y)} u_k (y) dy, \label{up} \\
\mbox{and }0 = \int_{\p \Omega^-} G_k(x,y) \f{\p u_k}{\p n} (y) 
- \f{\p  G_k(x,y)}{\p n(y)} u_k (y) dy, \label{down}
\eean
where $n$ is the exterior normal vector in each case. If $y$ is in $\p \Omega^+$ and is such that $y_2 = 0$ it is clear that $\ds \f{\p  G_k(x,y)}{\p n(y)} =0$. We also use that $u_k$ is even in $x_2$ so that 
$\ds \f{\p u_k}{\p x_2} (x)=0$ if $x_2 =0$ and $x \notin \Gamma$, since $u_k $ is even in $x_2$.
Combining (\ref{up}) and (\ref{down}) we find that for $x$ in $\Omega^+$
\bea
u_k (x) =  \int_{\p D} G_k(x,y) \f{\p u_k}{\p n} (y) 
- \f{\p  G_k(x,y)}{\p n(y)} u_k (y) dy \q
- \int_{\Gamma} G_k(x,y) ( \f{\p u_k}{\p x_2^+}  - \f{\p u_k}{\p x_2^-})
(y) dy
\eea
But due to point (iii) in lemma \ref{onuk}, since $u_k = \tilde{u}_k$ on $\p D$,
$$
 \int_{\p D} G_k(x,y) \f{\p u_k}{\p n} (y) 
- \f{\p  G_k(x,y)}{\p n(y)} u_k (y) dy \q =0,
$$
for all $x$ in $\Omega^+$, thus
\bean \label{ukintalmost}
u_k (x) = 
- \int_{\Gamma} G_k(x,y) ( \f{\p u_k}{\p x_2^+}  - \f{\p u_k}{\p x_2^-}) (y) dy,
\eean
for all $x$ in $\Omega^+$. 
Now take the $x_2$ derivative of (\ref{ukintalmost}) for $x$ in $\Omega^+$
approaching $\Gamma$. Due to the jump condition
for normal derivatives of single layer potentials,
we find that
\bea
 \f{\p u_k}{\p x_2^-} (x) = -\f12 ( \f{\p u_k}{\p x_2^+}  - \f{\p u_k}{\p x_2^-})(x)
- \int_{\Gamma} \f{\p G_k(x,y)}{\p x_2} ( \f{\p u_k}{\p x_2^+}  - \f{\p u_k}{\p x_2^-})
(y) dy.
\eea
We observe that for $x$ and $y$ on $\Gamma$, $\ds \f{\p G_k(x,y)}{\p x_2} = 0$. It follows that $\ds \f{\p u_k}{\p x_2^-} (x) = - \f{\p u_k}{\p x_2^+} (x)$ for $x$ on $\Gamma$.
$\square$

\begin{lemma}
$u_0$ is equal to the constant function 1.
\end{lemma}
\textbf{Proof:}
It can be shown, as in the case where $k>0$, that
 $u_0$ is in $H^1(\Omega)$,
the upper and lower trace on $\Gamma$ of $u_0$
are both equal to the constant 1, and
$u_0$ can be extended to a function in $\RR^2 \setminus\Gamma$
such that, if we still denote by $u_0$ the extension, 
$ \Delta  u_0=0$ in   $\RR^2 \setminus\Gamma$. 
The condition at infinity for $u_0$ is different. From 
(\ref{userieslap}) we infer that $u_0 = a_0 + O(r^{-1})$
and $ \p_r u_0 =  O(r^{-2})$.
We also note that (\ref{userieslap}) implies that 
\bean \label{zeroav}
 \int_{\p D}  \p_r {u}_0 =0
\eean
Since $\Delta u_0=0$ in $\Omega$, $u_0=1$ on $\Gamma$,
and (\ref{zeroav}) holds, we have
\bean
  \int_{\Gamma} \f{\p u_0}{\p x_2^+}  - \f{\p u_0}{\p x_2^-}  = 0
\eean
but the latter is also equal to
\bean
  \int_{\Gamma} (\f{\p u_0}{\p x_2^+}  - \f{\p u_0}{\p x_2^-}) 
\ov{u_0}  = 0
\eean
so applying Green's formula in combination to the estimates $u_0 = a_0 + O(r^{-1})$ and $ \p_r u_0 =  O(r^{-2})$ we find that $\ds \int_{\RR^2} |\nabla u_0|^2 =0$. We infer  $u_0$ is a constant in $ \RR^2$. That constant can only be 1.
$\square$

\begin{figure}
\begin{center}
\includegraphics[scale=.5]{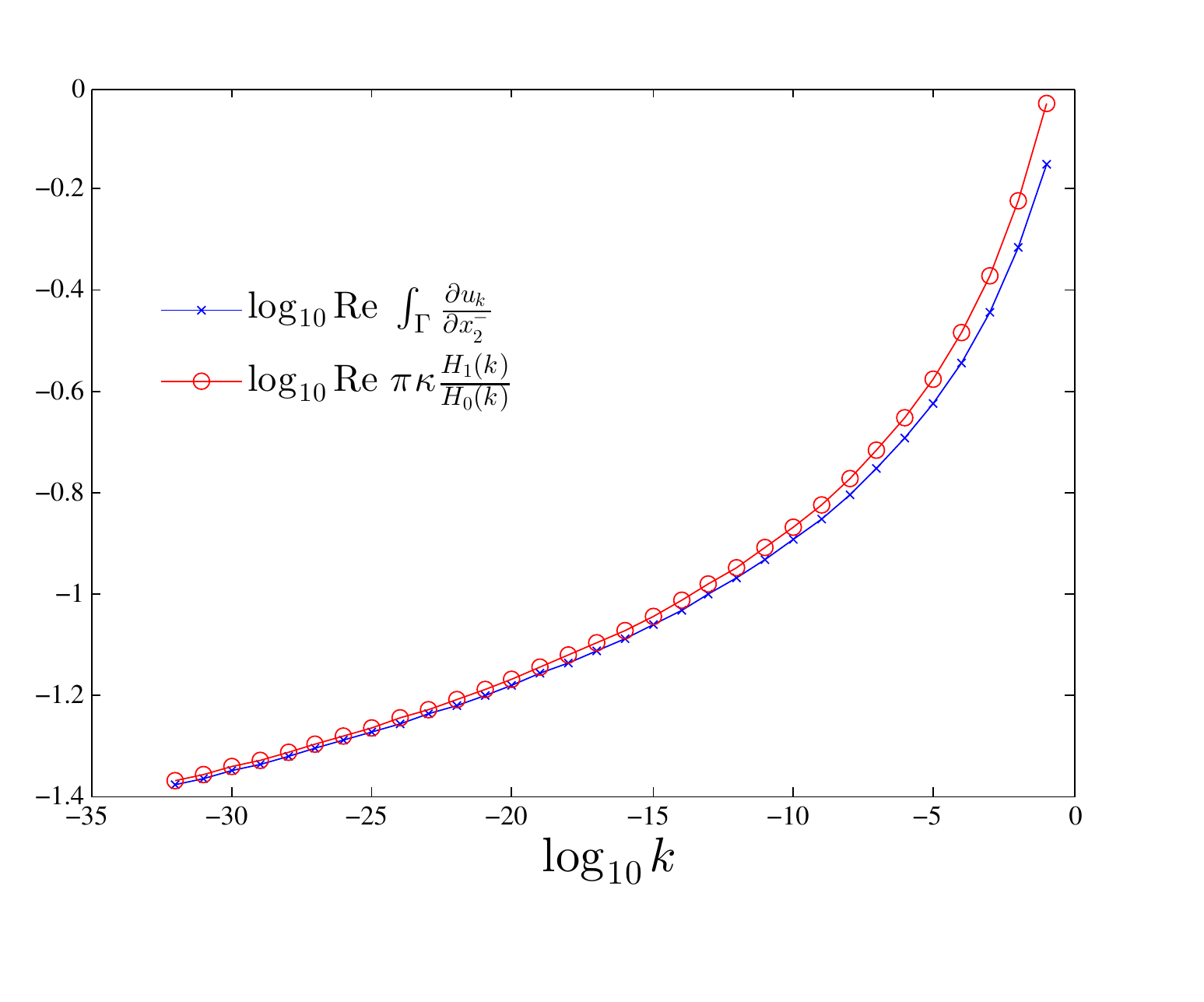}
\caption{$\ds \log_{10} \mbox{Re }\int_{\Gamma} \frac{\partial u_{k}}{\partial x_2^{-}}$ (blue) against $\ds \log_{10} 
 \mbox{Re }\big( \pi k \frac{H_1(k)}{H_0(k)} \big)$ (red). 
} 
\end{center}
\label{fig:small_wave_numbers_asympt}
\end{figure}

\begin{theorem} \label{low_k_equiv}
The following estimate as $k$ approaches $0^+$ holds
\bean \label{at0est}
\int_\Gamma \f{\p u_k}{\p x_2^-} \sim \pi k \f{H_1(k)}{H_0(k)}
\eean
Consequently, $\ds Re \int_\Gamma \f{\p u_k}{\p x_2^-}$ must be strictly positive for small values of $k > 0$.
\end{theorem}
\textbf{Proof:}
Set $v=1- \phi$  in variational problem (\ref{varphi})
(note that the trace of $v$ is zero on $\Gamma$, as required in the space $H^1_{0,\Gamma}(\Omega)$), to obtain
\bean
-\int_\Omega \nabla\tilde{u}_k \cdot \nabla \phi  - 
\int_\Omega  k^2 \tilde{u}_k (1 -\phi)-
\int_{\p D} (T_k \tilde{u}_k)\, (1 - \phi) = \int_\Omega (\Delta \phi + k^2\phi)(1- \phi),
\eean
We first observe that $\ds \int_\Omega k^2\phi (1- \phi) = O(k^2)$ and due to theorem \ref{convth}, $\ds \int_\Omega  k^2 \tilde{u}_k (1 -\phi)= O(k^2) $.\\
As $\phi$ is zero on $\p D$, we have found that
\bean \label{est1}
-\int_\Omega \nabla\tilde{u}_k \cdot \nabla \phi 
+ \int_\Omega \Delta \phi \,\phi = \int_{\p D} T_k \tilde{u}_k + O(k^2)
\eean
Next, using Green's theorem,
\bean \label{est2}
2 \int_\Gamma \f{\p \tilde{u}_k}{\p x_2^-} = 2 \int_\Gamma \f{\p \tilde{u}_k}{\p x_2^-} \phi =
\int_\Omega \nabla \tilde{u}_k \cdot  \nabla \phi  +
\int_\Omega  \Delta \tilde{u}_k \, \phi
\eean
Since in $\Omega$, $\Delta \tilde{u}_k = -\Delta \phi - k^2 \tilde{u}_k-k^2 \phi$ combining (\ref{est1}) and  (\ref{est2}) yields
\bean
2 \int_\Gamma \f{\p \tilde{u}_k}{\p x_2^-} = 
- \int_{\p D} T_k \tilde{u}_k + O(k^2)
\eean
But $\ds \int_{\p D} T_k \tilde{u}_k  = -k \f{H_1(k)}{H_0(k)} 2 \pi a_0(k)$
where $\ds a_0(k) =  \f{1}{2 \pi}\int_{\p D} \tilde{u}_k $,
so using again that $\tilde{u}_k  $ is strongly convergent to $1 - \phi$
in $H^1(\Omega)$, $a_0(k)$ tends to $1$ as $k \ri 0$, 
thus we claim that 
\bean
\int_{\p D} T_k \tilde{u}_k  \sim  -2 \pi k \f{H_1(k)}{H_0(k)} 
\eean
as $k \ri 0$. Going back to the definition of Hankel functions it easy to see that
 $\ds k \f{H_1(k)}{H_0(k)} \sim - (\ln k)^{-1}$ as $k \ri 0$,
from where we conclude that 
\bean \label{log_estimate}
\ds \int_\Gamma \f{ \p \tilde{u}_k}{\p x_2^-} \sim - \pi (\ln k)^{-1},
\eean 
so  $\ds \mbox{Re }\int_\Gamma \f{\p u_k}{\p x_2^-}$ must be strictly positive for all $k>0$ small enough. $\square$ \\
For illustration, in figure \ref{fig:small_wave_numbers_asympt} we have plotted graphs 
of $\ds \log_{10}\mbox{Re }\int_{\Gamma} \frac{\partial u_{k}}{\partial x_2^{-}}$ and 
$\ds \log_{10} \mbox{Re }\big( \pi k \frac{H_1(k)}{H_0(k)} \big)$ against  $\ds \log_{10} k$.
Note that
$\ds \int_{\Gamma} \frac{\partial u_{k}}{\partial x_2^{-}}$
was computed using the numerical method outlined in \cite{DVandSergey_computation} and \cite{ghergu}.

\section{Asymptotics for high wavenumbers}
We provide in this section a derivation of an equivalent for $\ds \int_\Gamma \f{\p u_k}{\p x_2^-}$ as $k \ri \infty$, where $u_k = \tilde{u}_k + \phi$, and $\tilde{u}_k$ solves variational problem (\ref{varphi}).
We will prove the following theorem
\begin{thm} \label{high_k_theorem}
Let $\tilde{u}_k$ be the solution to (\ref{varphi}),
and set $u_k = \tilde{u}_k + \phi$ . \\
The following estimates hold as $k$ approaches infinity,
\bean \label{main_est_k_large}
\int_\Gamma \f{\p u_k}{\p x_2^-} = -ik + O(k^{3/4}).
\eean
\end{thm}
What is the asymptotic behavior of $\f{\p u_k}{\p x_2^-}$ as the wavenumber $k$
approaches infinity? 
High frequency approximation for the wave equation is a vast subject which has been
extensively studied over time. Historically, investigators have tried to explain how the laws
of geometric optics relate to the wave equation at high frequency in an attempt to provide
a sound foundation for Fresnel's laws. Kirchhoff may have been the first one to write
specific equations and asymptotic formulas for  high frequency wave phenomena, however,
his derivation was informal. 
A more  mathematically rigorous  study of the behavior of solutions to the wave equation
at high frequency requires the use of Fourier integral operators and micro local analysis.
As far as we know this kind of work was pioneered by Majda, Melrose and Taylor, see 
\cite{Majda, MajdaTaylor, MelroseTaylor}. These authors were actually interested 
in the case of the exterior of a bounded convex domain, so their results can not be applied to our 
case since $\Gamma$ has empty interior in $\RR^2$. We have instead to rely on recent
groundbreaking work by 
Hewett, Langdony, and, Chandler-Wildez, see \cite{Chandler1, Chandler2} which pertains 
to either scattering in dimension 2 by soft or hard line segments (our case), or
scattering in dimension 3 by soft or hard  open planar surfaces.
The great achievement of their work is that they were able to derive continuity
and coercivity bounds that explicitly depend  on the wavenumber.\\
Following the work by Hewett and Chandler-Wildez
we introduce relevant functional spaces and frequency depending norms.
Let $v$ be a tempered distribution on $\RR$ and $\hat{v}$ its Fourier transform. 
Let $s$ be in $\RR$.  We say that $v$ is in $H^s(\RR)$ if
$\int_\RR (1 + \xi^2)^s |\hat{v} (\xi)|^2 d \xi< \infty$.
We then define in $H^s(\RR)$  the $k$ dependent norm
\bean
\| v \|_{H^s_k(\RR)} = (\int_\RR (k^2 + \xi^2)^s |\hat{v} (\xi)|^2 d \xi)^{1/2}.
\eean
Note that $H^1(\RR)$ is included in $H^{1/2}(\RR)$ and more precisely,
\bean
\label{Sobolev_inj}
 \| v \|_{H^{1/2}_k(\RR)} \leq k^{-1/2} \| v \|_{H^{1}_k(\RR)},
\eean
for all $v$ in $H^1(\RR)$.
Denote by $I$ the interval $\ds(-\f12, \f12)$.
$\tilde{H}^s(I)$ is defined to be the closure of $C^\infty_c(I)$ 
(which is the space of smooth functions, compactly supported in $I$)
for the norm
$\|  \|_{H^s_k(\RR)}$. 
${H}^s(I)$ is defined to be the space of restrictions to $I$ of elements
in $H^s(\RR)$. We define on ${H}^s(I)$ the norm
\bea
 \| v \|_{H^s_k(I)} = \inf \{ \| V\|_{H^s_k(\RR)}: V \in H^s_k(\RR) \mbox{ and }
V|_I = v \}
\eea
\begin{thm} \label{Chandler} (Hewett and Chandler-Wildez)
For any $s$ in $\RR$, the
 operator $S_k$ defined by the following formula for smooth functions $v$ on $I$
$$
(S_k v) (x_1) = \int_I \f{i}{4} H_0 (k |x_1 - y_1|) v (y_1)  dy_1
$$
can be extended to a continuous linear operator
from $\tilde{H}^s(I)$ to $H^{s+1}(I)$. 
Furthermore, $S_k$ is injective, $S_k^{-1}$ is continuous and satisfies the estimate
\bean \label{main_estimate}
 \| S_k^{-1}  v \|_{\tilde{H}^{-1/2}_k(I)} \leq 2 \sqrt{2} \|  v \|_{{H}^{1/2}_k(I)},
\eean
for all $v$ in $H^{1/2}_k(I)$.
\end{thm}
Let us emphasize one more time 
that although the continuity and coercivity properties of $S_k$ have been known for some
time, Hewett and Chandler-Wildez's great achievement was to derive the dependency of the 
   coercivity and the continuity bounds on the wavenumber $k$ as in estimate 
	(\ref{main_estimate}): the dependency appears in the use of the special norms
	$\|  . \|_{{H}^{s}_k(I)}$.

	\subsection{An informal derivation of estimate (\ref{main_est_k_large})}
	This informal derivation will be helpful since it will give us an idea of 
	what the asymptotic behavior of $\int_\Gamma \f{\p u_k}{\p x_2^-}$ should be.
	 Denote by $f_k(y_1)$ the value of $\ds \f{\p u_k}{\p x_2^-} (y_1, 0)$ for $y_1$ in $I $.
	According to (\ref{intform}), $f_k $ satisfies the integral equation $S_k f_k = \f12$ or
\bean \label{inteq} 
\int_I \f{i}{4} H_0 (k |x_1 - y_1|) f_k (y_1) d y_1 = \f12, \q
x_1 \in I
\eean
We multiply equation (\ref{inteq}) by $-2ik$ and we integrate in $x_1$ over $I$ to obtain
\bean \label{integrated} 
\int_I \int_I \f{k}{2}  H_0 (k |x_1 - y_1|) f_k (y_1) d y_1 d x_1 = - ik
\eean
If it is possible to interchange the order of integration in the left hand side of
(\ref{integrated}) then
\bean \label{integrated2} 
\int_I \int_I \f{k}{2} H_0 (k |x_1 - y_1|) d x_1 f_k (y_1) d y_1 = - ik
\eean
Now for every $x_1$ in $I$, 
 setting $v=k(x_1 - y_1)$ and  using lemma \ref{struve} from the appendix, we find that
$$
\ds \lim_{k \ri \infty} \int_I \f{k}{2} H_0 (k |x_1 - y_1|) d x_1 =1,
$$
so we are led to believe that $ \int_I  f_k (y_1) d y_1 \sim - ik$, which we  set out
to prove in the next section.

\subsection{Rigorous derivation of estimate (\ref{main_est_k_large})}
\begin{lemma}	\label{firstHknorm}
	The following estimate holds as $k$ approaches infinity
\bean \label{estimateik} 
	\| S_k(-ik) -\f12 \|_{{H}^{1/2}_k(I)} = O(k^{1/4})
\eean
\end{lemma}	
\textbf{Proof:}
We start by recalling lemma \ref{struve} from the appendix and noting that for $x_1$ in I,
$$
S_k(-ik)(x_1) -\f12 = \int_I \f{1}{4} H_0 (k |x_1 - y_1|) k  dy_1 - \f12
= -\int_{k(\f12 + x_1)}^\infty \f14 H_0 (v) dv -\int_{k(\f12 - x_1)}^\infty \f14 H_0 (v) dv
$$
We now set $g_k(x_1) = \int_{k(\f12 + x_1)}^\infty H_0 (v) dv + \int_{k(\f12 - x_1)}^\infty  H_0 (v) dv$.
Since the function $t \ri \int_{t}^\infty H_0 (v) dv$ is continuous on $[0, \infty)$ and has limit 
zero at infinity, there is a positive $C$ such that 
\bean \label{est1forH}
|\int_{t}^\infty H_0 (v) dv| \leq C, 
\eean
for all $t$ in $[0, \infty)$.
It is well known that as $v$ approaches 
infinity (see  \cite{Abra}), $ H_0 (v) = e^{i(v- \f{\pi}{4})} \sqrt{\f{2}{\pi v}} + O(v^{-\f32})$,
so we also have the estimate
\bean \label{est2forH}
|\int_{t}^\infty H_0 (v) dv| = O(t^{-\f12}),  \q t \ri \infty
\eean
Without loss of generality we may assume that $k>4$.
 If $x_1$ is in $(-\f12 + k^{-1/2}, \f12 - k^{1/2})$,
then  $k(\f12 + x_1)$ and $k(\f12 - x_1)$ are greater than 
$k^{1/2}$ so due to
(\ref{est2forH}), $g_k(x_1) = O(k^{-1/4})$. 
Then using (\ref{est1forH}) we infer that 
\bean \label{L2conv}
\int_I |g_k|^2 =  \int_{(-\f12 + k^{-1/2}, \f12 - k^{1/2})} |g_k|^2
+ \int_{I \setminus (-\f12 + k^{-1/2}, \f12 - k^{1/2})} |g_k|^2 = O(k^{-1/2})
\eean
Next we note that $g_k'(x_1) = k H_0(k(\f12 + x_1)) - k  H_0(k(\f12 - x_1))$.
Using a substitution we find that
$$
\int_I |g_k'|^2 \leq  4k  \int_0^k |H_0(v)|^2 dv
$$
Given the asymptotics at infinity of $H_0$ we infer
that 
\bean
\label{L2priconv}
\int_I |g_k'|^2  = O( k \ln k ), \q k \ri \infty 
\eean
We now use (\ref{L2conv}) and (\ref{L2priconv}) to evaluate a few 
Sobolev norms of $g_k$. 
We observe that $g_k(-\f12)=g_k(\f12)=\int_{0}^\infty H_0 (v) dv + \int_{k}^\infty  H_0 (v) dv$
is uniformly bounded in $k$. 
We set for $k>0$
$$
\tilde{g}_k(x_1) = \left\{
\begin{array}{ll}
0 & \mbox{ if } |x_1| \geq \f12 + k^{-1/2} \\
k^{1/2} (x +\f12  + k^{-1/2} ) g_k(\f12) &\mbox{ if }  -\f12 - k^{-1/2}<x_1 \leq - \f12 \\
g_k(x_1) &\mbox{ if } -\f12 < x_1 \leq \f12 \\
-k^{1/2} (x -\f12  - k^{-1/2}) g_k(\f12) & \mbox{ if }  \f12<x_1 <\f12 + k^{-1/2}
\end{array}
\right.
$$
Clearly, $\tilde{g}_k$ is in $H^1(\RR)$.
A straightforward calculation will show that 
\bean
\label{tilde_est}
\int_{\RR \setminus I} |\tilde{g}_k|^2 = O(k^{-1/2})
\mbox{ and } \int_{\RR \setminus I} |\tilde{g}_k'|^2 = O(k^{1/2})
\eean
Estimates (\ref{tilde_est}) 
combined to (\ref{L2conv}) and (\ref{L2priconv})
implies the following estimates for the Fourier transform of $\tilde{g}_k$
\bean
\label{fourier_est}
\int_{\RR } |\widehat{\tilde{g}_k}(\xi)|^2 d \xi = O(k^{-1/2})
\mbox{ and } \int_{\RR } \xi^2|\widehat{\tilde{g}_k}(\xi)|^2 d \xi  = O(k \ln k),
\eean
so 
$$
\int_{\RR } (k^2+\xi^2)|\widehat{\tilde{g}_k}(\xi)|^2 d \xi  = O(k^{3/2}),
$$
that is, $\| \tilde{g}_k \|_{{H}^{1}_k(\RR)} = O(k^{3/4})$.
Now,  due to inequality (\ref{Sobolev_inj}), $\| \tilde{g}_k \|_{{H}^{1/2}_k(\RR)} = O(k^{1/4})$
and estimate (\ref{estimateik}) is proved. $\square$

\begin{lemma}	\label{secondHknorm}
	The following estimate holds as $k$ approaches infinity
\bean \label{estimateof1} 
	\| 1 \|_{{H}^{1/2}_k(I)} = O(k^{1/2} )
\eean
\end{lemma}	
\textbf{Proof:}
We set for $k>0$
$$
h(x_1) = \left\{
\begin{array}{ll}
0 & \mbox{ if } |x_1| \geq 1\\
2(x +1 )  &\mbox{ if }  -1<x_1 \leq - \f12 \\
1 &\mbox{ if } -\f12 < x_1 \leq \f12 \\
-2 (x -1)  & \mbox{ if }  \f12<x_1 <1
\end{array}
\right.
$$
Clearly $h$ is in $H^1(\RR)$ and is independent of $k$ so
\bean
\label{fourier_est3}
\int_{\RR } |\widehat{h}(\xi)|^2 d \xi = O(1)
\mbox{ and } \int_{\RR } \xi^2|\widehat{h}(\xi)|^2 d \xi  = O(1)
\eean
From there we infer  $\| h \|_{H^{1/2}_k(\RR)} = O(k^{1/2})$
and estimate (\ref{estimateof1}) is proved. $\square$ \\\\
\textbf{Proof of theorem \ref{high_k_theorem}:}\\
Combining estimates (\ref{main_estimate}) and (\ref{estimateik})
and recalling that $S_k f_k = \f{1}{2}$ we may write
\bean
\| -ik - f_k \|_{\tilde{H}^{-1/2}_k(I)} = O(k^{1/4})
\eean
We then use estimate (\ref{estimateof1}) and we write denoting by $<,>$ the duality 
bracket $\tilde{H}^{-1/2}_k(I), H^{1/2}_k(I)$ 
\bean \label{O34}
<-ik -f_k, 1> = O(k^{3/4})
\eean
Following \cite{Hsiao_Stephan_Wendland}, there are two constants 
$C_1$ and $C_2$
such that 
$$f_k - C_1 (\f12 - x_1)^{-1} - C_2 (\f12 + x_1)^{-1} 
$$
is in $H^{1/2}_k(I)$.
In particular $f_k$ is in $L^q(I)$ for every $q$
in $[1,2)$, so $\int_I f_k$ is a regular Lebesgue integral.
Since $I$ has measure 1,
estimate (\ref{main_est_k_large}) follows from (\ref{O34}).  $\square$ \\\\

\begin{figure}
\begin{center}
  \includegraphics[scale=.4]{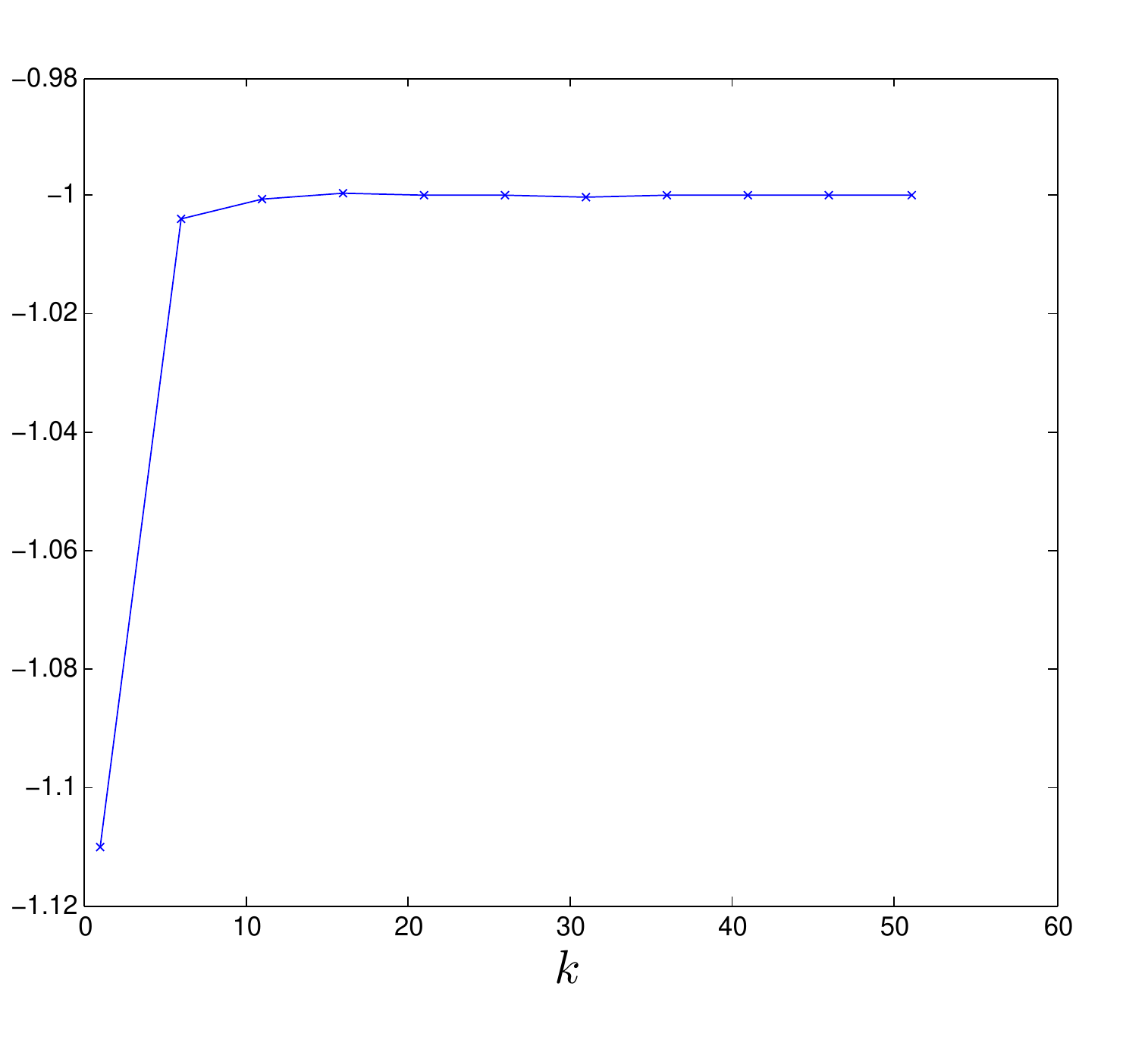} \includegraphics[scale=.4]{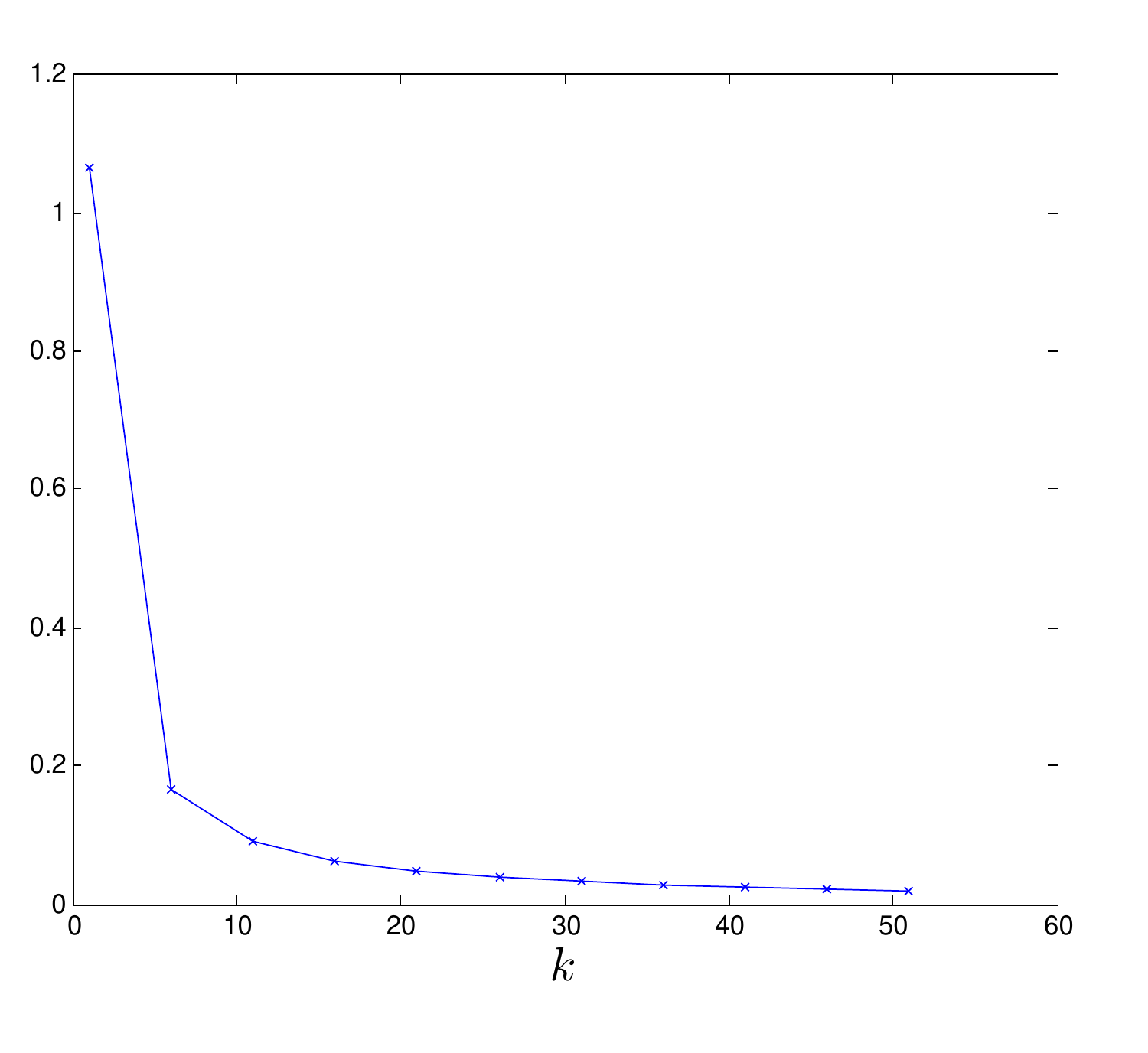}
	\end{center}
  \caption{Left: computed values of  $\ds \frac{1}{k} Im \int_\Gamma \f{\p u_k}{\p x_2^-}  $ 
	against $k$. Right: computed values of $\ds \frac{1}{k} Re \int_\Gamma \f{\p u_k}{\p x_2^-}$
	against $k$.
	}
  \label{fig:high_wave_numbers_asympt}
\end{figure}

For illustration, figure \ref{fig:high_wave_numbers_asympt} shows a graph of the imaginary part 
and a graph of the real part of
$\ds \frac{1}{k} \int_{-\frac{1}{2}}^{\frac{1}{2}} f_{k} (y_1) dy_1$,
that is, $\ds \frac{1}{k} \int_\Gamma \f{\p u_k}{\p x_2^-}$,
 against the wavenumber $k$. 
See \cite{DVandSergey_computation} and \cite{ghergu}
for the numerical method that we used.

\section{Proof of theorem \ref{main_th}}

\subsection{The case $C_1, C_2, C_3, C_4 >0$}
\begin{figure}[htb]
\begin{center}
  \includegraphics[scale=.4]{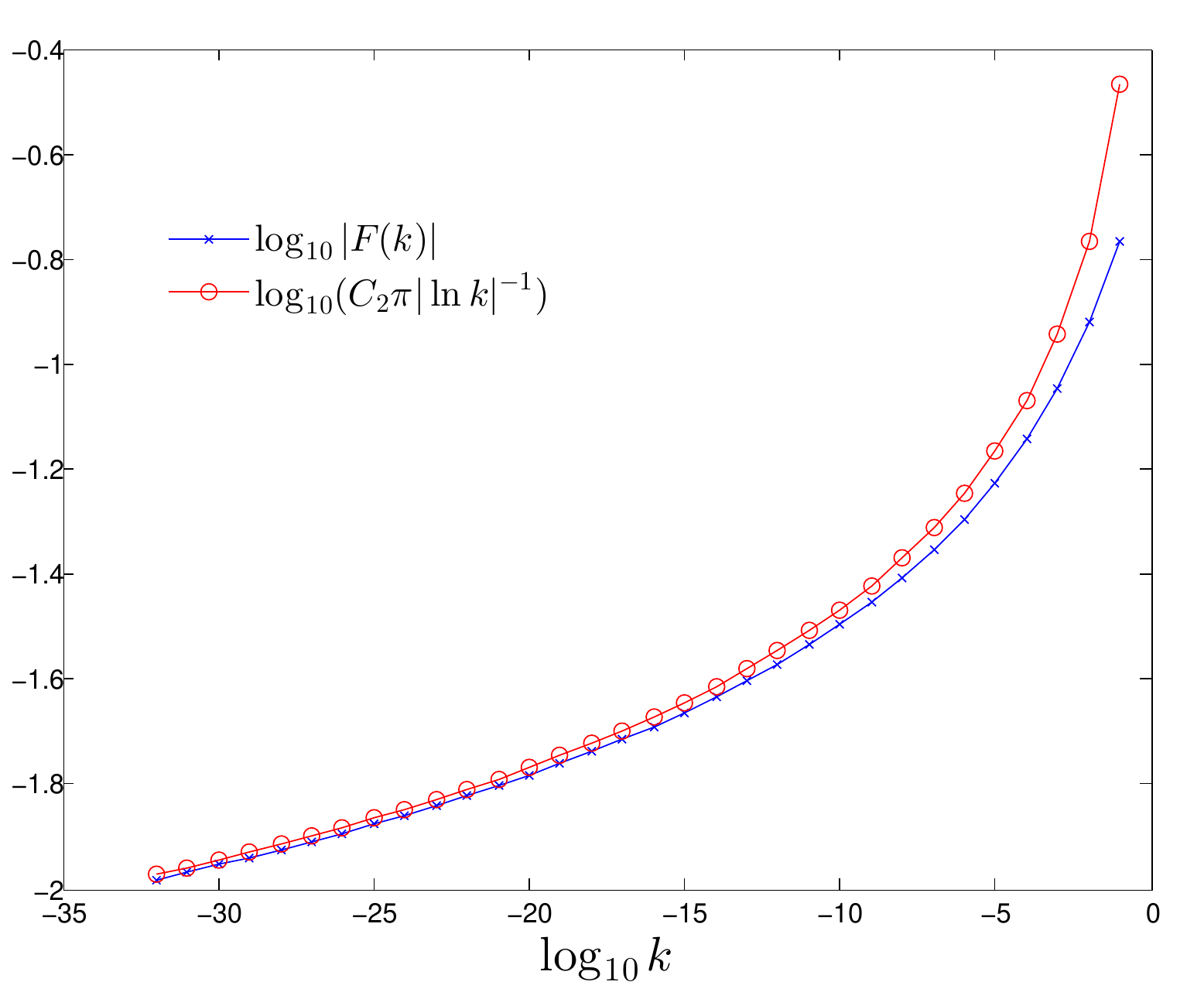} \includegraphics[scale=.4]{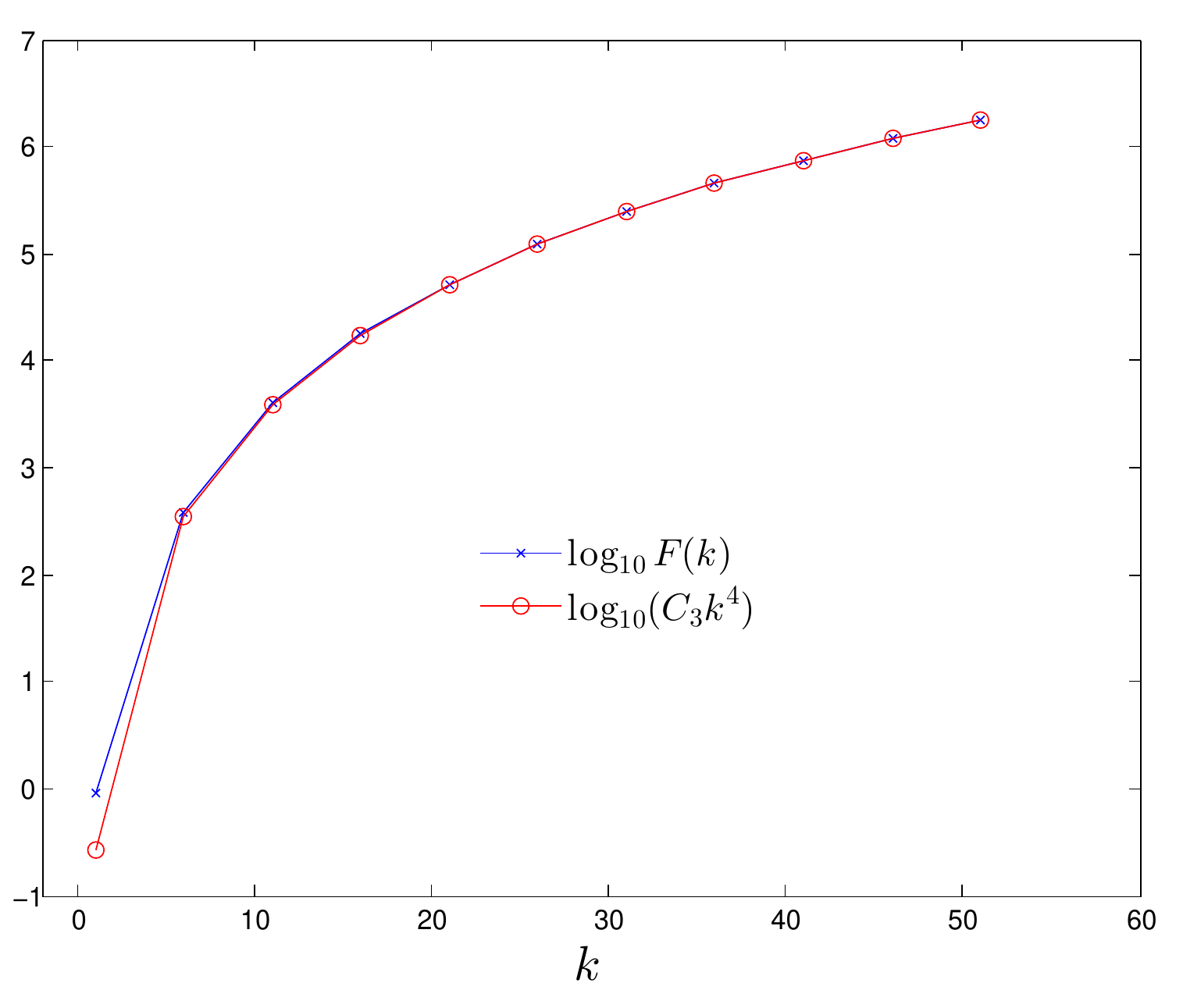}
	\end{center}
  \caption{ Numerical computations illustrating formulas (\ref{atzero}) and (\ref{atinfinity}).
	These graphs were produced using the specific constants  $C_1 = 0.4, C_2 = 2/3, C_3 = 5/12, C_4 = 5/48 $,
see \cite{Sergey_sthesis} and \cite{DVandSergey_computation} for why this choice is physically plausible and for 
a description of the numerical method employed.
	Left: $\log_{10} |F(k)|$ (blue) vs $\log_{10} (C_2 \pi |\ln k|^{-1})$ (red). Right: $\log_{10} F(k) $ 
	(blue) vs $\log_{10} (C_3 k^4)$ (red).}   
  \label{fig:p_q_asympt}
\end{figure}

In this section we cover the case where the constants $C_1, C_2, C_3, C_4 $ are all positive, which was assumed 
in theorem \ref{main_th}.
We explained in section \ref{section_equivalent_problem}
how the function $\Phi$ defined by 
(\ref{non_dim_helmholtz}-\ref{uext0}) relates to $u_k$ where 
 $u_k = \tilde{u}_k + \phi$ and $\tilde{u}_k$ solves (\ref{varphi}):
if $u_k$ is extended to $\RR^2$ as indicated by lemma \ref{onuk} then
$u_k$ and $\Phi$ are equal on $\ov{\RR^{2+}}$.
Accordingly, recalling   the definition of the function $F$ given in (\ref{Fdef}) and formula (\ref{opp}),
$F$ can also be expressed as follows,
\bean
F(k)=q(k^2) + p(k^2)\mbox{Re }\int_\Gamma \f{\p u_k}{\p x_2^-}
\eean
We know  from lemma \ref{convth}
  that  $F$ is analytic in $(0, \infty)$.
Recalling the definition of $p$ and $q$, 
(\ref{p_and_q}),
and using  estimate (\ref{log_estimate}) we have that 
\bean \label{atzero}
F(k) \sim  C_2 \pi  (\ln k)^{-1}, \q k \ri 0^+
\eean
Recalling estimate (\ref{main_est_k_large}) we claim that
\bean  \label{atinfinity}
F(k) \sim  C_3 k^4, \q k \ri \infty
\eean
It follows from estimate (\ref{atzero}) that there is a positive
$\alpha$ such that $F(k)<0$ if $k$ is in $(0,\alpha)$,
and  from estimate (\ref{atinfinity}) we infer that
$\lim_{k \ri \infty} F(k) =  \infty$. Since $F$ is continuous 
in $(0,\infty)$,  we conclude that $F$ must achieve the value zero on that interval. We can also claim that the zeros of $F$ are isolated since $F$ is an analytic function. Due to (\ref{atinfinity}) these zeros occur only in some interval $[A,B]$ where $A$ and $B$ are two positive constants. In particular the equation $F(k)=0$ has at most a finite number of solutions. \\
For illustration, figure \ref{fig:p_q_asympt} shows $\log_{10} |F(k)|$ and $\log_{10} (C_2 \pi |\ln k|^{-1})$ against $\log_{10} k$ on the left, and $\log_{10} \big| F(k) \big|$ and $\log_{10} (C_3 k^4)$ against $k$ on the right. These graphs were produced using the specific constants $C_1 = 0.4, C_2 = 2/3, C_3 = 5/12, C_4 = 5/48 $,
see \cite{Sergey_sthesis} and \cite{DVandSergey_computation} for why this choice is physically plausible.

\subsection{Use of estimates in the case where $C_3<0$}
We first claim that for any positive value of $C_1, C_2, C_4$, \textsl{there exist  negative values
of $C_3$, such that the equation $F(k)=0$ has no solutions}.
According to estimate (\ref{log_estimate}) there is a positive $\alpha$ such that for all $k$ in 
$(0, \alpha)$,
\bean \label{from0toalpha}
- \f12 \pi (\ln k)^{-1} \leq 
 \mbox{Re }\int_\Gamma \f{\p u_k}{\p x_2^-} \leq - \f32\pi (\ln k)^{-1}.
\eean
Consequently for $k$ in $(0, \alpha)$, 
\bean
F(k) \leq C_4 k^2 + C_1 k^2 (- \f32\pi (\ln k)^{-1}) + C_2 \f12 \pi (\ln k)^{-1}  
\eean
so if $\alpha$ is small enough, $F(k)<0$ if $k$ is in $(0, \alpha)$.
 According to (\ref{main_est_k_large}) there is a negative $A_1$ and a positive $A_2$
such that for all k in $[\alpha, \infty)$,
\bean \label{past_alpha}
A_1 k^{3/4} \leq 
 \mbox{Re } \int_\Gamma \f{\p u_k}{\p x_2^-} \leq A_2 k^{3/4}.
\eean
so 
 for all k in $[\alpha, \infty)$,
\bean \label{past_alpha2}
F(k) \leq C_3 k^4 + C_4 k^2 + C_1 A_2 k^{2+3/4} - C_2 A_1 k^{3/4}
\eean
so if we choose $C_3$ less than some negative constant, $F(k)<0$ for all k in $[\alpha, \infty)$.
We conclude that for that choice of $C_3$, $F(k)<0$ for all positive $k$,
so the equation $F(k)=0$ has no solution in $(0, \infty]$. \\

Next we show the following claim: for any value of the positive constants $C_1, C_2$ and any value of the negative constant 
$C_3$, 
\textsl{for all values of $C_4$ greater than some constant, the equation $F(k)=0$ has at least one solution}.
We may assume that $\alpha$ defined above is less than 1.
For  all k in $[\alpha, \infty)$
\bean
F(k) \geq  C_3 k^4 + C_4 k^2  + C_1 A_1 k^{2+3/4} - C_2 A_2 k^{3/4}
\eean
thus $F(1)>0$ for any $C_4$ greater than some constant. As $C_3 <0$, 
thanks to (\ref{past_alpha2}), we have that,
$\lim_{k \ri \infty} F(k) =  - \infty$, so by continuity of $F$ we conclude that the equation
$F(k)=0$ has at least one solution in $(0, \infty)$.

\section{Conclusion and perspectives}
We have proved in this paper the existence of
frequency modes  coupling  seismic waves and vibrating tall buildings. 
 Although the derivation from physical principles of a set of equations  modeling this phenomenon
was carefully done and numerical evidence supporting existence of coupling
frequency modes was shown in previous studies, to the best of our knowledge, 
 up to now,
 there was  no  formal proof of existence. \\
We believe that only minor modifications to our present argument would show existence 
of coupling frequency modes in the case of a finite set of buildings. 
The present study was limited to a single building chiefly for clarity of exposition.
It could turn out, however, that the case of a periodically repeated pattern of buildings might involve
more substantial alterations of the existence proof given in this paper. Note that this case was particularly important
in our simulations of coupling modes presented in \cite{DVandSergey_computation}. \\
Recalling that our study pertains exclusively to anti plane shearing, it will be important to generalize 
our results to fully three dimensional elastic vibrations. At present, this appears to be quite an undertaking given
the complexity of  Green's tensor for half space  elasticity and the lack
of a known analog to theorem \ref{Chandler} for the elasticity operator.

\section{Appendix: useful results on Hankel functions}
For $n \in \NN$ and $z$ in $\CC$ we denote by $J_n$  the Bessel function
$$
J_n(z) = (\f12 z)^n \sum_{k=0}^\infty (-\f14)^k \f{z^{2k}}{k! (n+k)!}
$$
For $n \in \NN$ and $z$ in $\CC\setminus \RR^{-*}$
we denote by $Y_n$ the Bessel function
\bea
-\f{1}{\pi} (\f12 z)^{-n}\sum_{k=0}^{n-1}
\f{(n-k-1)!}{k!} (\f14 z^2)^k +
\f{2}{\pi} \ln (\f12 z) J_n(z) \\
- \f{1}{\pi} (\f12 z)^{n}\sum_{k=0}^{\infty}
(\psi(k+1) + \psi(n+k+1) )(-\f14)^k \f{z^{2k}}{k! (n+k)!},
\eea
where the first sum is void if $n=0$ and the function $\psi$
is defined by $$\psi(1)=-\gamma, \,
\psi(n)= - \gamma + \sum_{k=1}^{n-1} \f{1}{k},$$ and $\gamma$
is the Euler constant.
 The Hankel function of the first kind of order $n$, that is,
$ J_n +i Y_n$, will be denoted
$H_n$.
\begin{lemma}\label{realeq}
The following equivalence as $n \ri \infty$ is uniform 
for all $z$ in a compact set of $(0,\infty)$
\bean
H_n(z) \sim - \f{i}{\pi} (\f12 z)^{-n} (n-1)!
\label{largen}
\eean
\end{lemma}
\textbf{Proof:}
Since for any positive $z$
\bea
|J_n(z) n! (\f12 z)^{-n} -1  | = \f{1}{n+1} |\sum_{k=1}^\infty (-\f14)^k \f{z^{2k}(n+1)!}{k! (n+k)!}| \\
\leq  \f{1}{n+1} \sum_{k=1}^\infty (\f14)^k \f{z^{2k}}{k!}
\eea
it is clear that $\ds J_n(z) \sim (\f12 z)^{n} \f{1}{n!}$ uniformly 
for all $z$ in a compact set of $(0,\infty)$. \\
We can show similarly that $\ds Y_n(z) \sim -\f{1}{\pi} (\f12 z)^{-n} (n-1)!$ uniformly 
for all $z$ in a compact set of $(0,\infty)$. 
We conclude that (\ref{largen}) must hold.
 $\square$

\begin{lemma}\label{realeq2}
For $z>0$, the following limit as $z \ri 0$ is uniform 
for all integers $n$ different from 0
\bean
\lim_{z \ri 0^+}-\f{ z}{|n|} \f{H_n'(z)}{H_n(z)} = 1
\label{ktozero}
\eean
For the special case $n=0$ we have,
$$
\lim_{z \ri 0^+}  \f{z H_0'(z)}{H_0(z)} = 0
$$
\end{lemma}
\textbf{Proof:}
We observe that for $n \geq 2$
$$\left( H_n(z) (- \f{i}{\pi} (\f12 z)^{-n} (n-1)!)^{-1}
-1 \right)
(n-1)$$
can be bounded by a function in $z$ which is continuous 
on $[0, \infty)$ and independent of $n$, therefore
\bean
H_n(z) \sim - \f{i}{\pi} (\f12 z)^{-n} (n-1)!
\label{largen2}
\eean
as $z \ri 0^+$, uniformly on $[0,b]$ for a fixed $b>0$.
Using the formula
\bean \label{besselrec}
H_n'(z)= - H_{n+1} (z) + \f{n}{z} H_n(z)
\eean
we obtain
$$
-\f{z}{n} \f{H_n'(z)}{H_n(z)} \sim 1
$$
as $z \ri 0^+$, uniformly on $[0,b]$ for a fixed $b>0$. \\
For $n \leq -2$, we can now apply the formula $H_{-n}(z) =(-1)^n H_n(z)$.
The remaining three cases can be treated in a straightforward fashion.
$\square$

\begin{lemma} \label{decrease}
For any $n$ in $\NN$, $|H_n(z)|$ is a decreasing function of $z$ on $(0,+\infty)$.
\end{lemma}
\textbf{Proof:}
This is due to the formula derived by Nicholson, see \cite{Watson},
$$
J_n^2(z) + Y_n^2(z) =
\f{8}{\pi^2} \int_0^\infty K_0(2 z \sinh t) \cosh 2nt dt
$$
where $\ds K_0(s) = \int_0^\infty e^{-s \cosh t} dt$.
 $\square$
\begin{lemma}\label{complexeq}
Let $a$ and $b$ be two real numbers such that $0<a<b$.
The following equivalence as $n \ri \infty$ is uniform 
for all  complex numbers $z$ in the closed disk
of the complex plane centered at $b$ and of radius $a$ 
\bean
H_n(z) \sim - \f{i}{\pi} (\f12 z)^{-n} (n-1)!
\label{largen3}
\eean
\end{lemma}
\textbf{Proof:}
The proof is nearly identical to that 
of lemma \ref{realeq}. $\square$

\begin{lemma}
\label{struve}
Denote by $\mbox{St}_n$ the Struve function of order $n$ as defined
in \cite{Abra}.
The following formula holds for any $t>0$
\bean
\int_0^t H_0(z) d z =t H_0(t) + \f{\pi}{2} t (\mbox{St}_0 H_1(t) - \mbox{St}_1 H_0(t))
\label{formst}
 \eean
It follows that the semi convergent integral $\ds \int_0^\infty  H_0(z) dz$ is exactly equal to $1$.
\end{lemma}
\textbf{Proof:}
Integral formula (\ref{formst}) is given in \cite{Abra}. 
The value of $\int_0^\infty  H_0(z) d z $ results from that formula combined to 
known asymptotics at infinity of Bessel and of Struve functions.
One should consult \cite{Abra} for formulas on Bessel functions,
and \cite{Watson} for their derivation.
$\square$


\begin{thebibliography}{99}
\bibitem{Abra} M. Abramowitz and I. Stegun, eds. (1992), Handbook of Mathematical Functions with Formulas,
Graphs, and Mathematical Tables, Dover, New York.



\bibitem{ullevi}
  S. Erlingsson, A. Bodare, Live load induced vibrations in Ullevi Stadium - dynamic underground analysis, \emph{underground Dyn. and Earthquake Eng.}, vol. 15, Issue 3 (1996) 171-188.


\bibitem{Folland}
G. B. Folland, 
Introduction to Partial Differential Equations, Second Edition, 
Princeton University Press.

\bibitem{ghergu}
  M. Ghergu, I. R. Ionescu, Structure-underground-structure coupling in seismic excitation and "city-effect". \emph{Int. J. Eng. Sci.} 47 (2009) 342-354.   


\bibitem{Chandler1}
D. P. Hewett,  S. N. Chandler-Wilde,
Acoustic scattering by fractal screens: mathematical formulations and wavenumber-explicit continuity and coercivity estimates,
arXiv.org > math > arXiv:1401.2805, 2014.

\bibitem{Chandler2}
 D. P. Hewett, S. Langdon, S. N. Chandler-Wilde,
A frequency-independent boundary element method for scattering by two-dimensional screens and apertures,
arXiv.org > math > arXiv:1401.2786, 2014.

\bibitem{Hsiao_Stephan_Wendland}
G. Hsiao, E. P. Stephan, W. L. Wendland,
On the Dirichlet problem in elasticity for a domain exterior to an arc,
Journal of Computational and Applied Mathematics, Volume 34, Issue 1, 10 February 1991, pp. 1-19

\bibitem{shuttles}
  H. Kanamori, J. Mori, B. Sturtevant, D. L. Anderson, T. Heaton, Seismic excitation by space shuttles, \emph{Shock waves}, 2 (1992) 89-96.


\bibitem{Kato}
T. Kato,  Perturbation theory for linear operators, Grundlehren Math. Wiss., Vol. 132. Springer, 1995.


\bibitem{Majda}
A. Majda, 
High frequency asymptotics for the scattering matrix and the inverse problem of acoustical scattering,
Communications on Pure and Applied Mathematics,
Volume 29, Issue 3, 261-291,  1976.

\bibitem{MajdaTaylor}
A. Majda,
M. E. Taylor,
The asymptotic behavior of the diffraction peak in classical scattering,
Communications on Pure and Applied Mathematics,
Volume 30, Issue 5,  639-669,  1977

\bibitem{MelroseTaylor}
Richard B Melrose, Michael E Taylor,
Near peak scattering and the corrected Kirchhoff approximation for a convex obstacle,
Advances in Mathematics 
Volume 55-3, 1985,  242-315.




\bibitem{Troianiello} Giovanni Maria Troianiello,
Elliptic Differential Equations and Obstacle Problems, University Series in Mathematics, 1987.

\bibitem{DVandSergey_computation} D. Volkov 
and S. Zheltukhin, Preferred Frequencies for Coupling of Seismic Waves and Vibrating Tall Buildings,
 submitted. (submission available on arXiv.org)


\bibitem{Watson} G. N. Watson, A treatise on the theory of Bessel functions, 
Cambridge University press, 1922.

\bibitem{Sergey_sthesis} S. Zheltukhin, Preferred Frequencies for Coupling of Seismic Waves and Vibrating Tall Buildings, PhD. thesis, WPI, 2013.

\end{thebibliography}
\end{document}